\input psfig.sty 
\magnification=\magstep1
\hsize=6.25truein
\overfullrule=0pt

\font\title=cmr10 scaled\magstep4
\font\subtitle=cmr10 scaled\magstep2

\font\author =cmr10 scaled\magstep1
\font\eaddr=cmtt10
\font\ninerm=cmr8
\baselineskip=14pt plus 1pt minus 1pt
\def\QED{\vrule height8pt width4pt depth0pt}

\def\proclaim #1: #2\par{\medbreak
	      \noindent{\bf#1:\enspace}{\sl#2}\par\medbreak}
\pageno=1	      

\hfill 21/9/9\break

\vskip 1truein
\centerline{\title Hard Core entropy:}
\vskip .1truein
\centerline{\title lower bounds}
\vskip .7truein
\centerline{\author Kari Eloranta}

\vskip .1truein
\centerline{\author Institute of Mathematics}
\centerline{\author Helsinki University of Technology}
\centerline{\author 02015 Espoo, Finland}
\vskip .1truein
\centerline{\eaddr kve@math.hut.fi}

\vskip .3truein

\centerline{\subtitle Abstract}
\vskip .3truein
\centerline{\vbox{\hsize 4.5in \noindent \ninerm \strut We establish
lower bounds for the entropy of the Hard Core Model on a few
2d lattices $\scriptstyle {\rm {\bf L}}.$ In this model the allowed
configurations inside $\scriptstyle \{0,1\}^{{\rm {\bf L}}}$ are the
one's in which the nearest neighbor $\scriptstyle 1$'s are forbidden.
Our method which is based on a sequential fill-in scheme is unbiassed and
thereby yields in principle arbitrarily good estimates for the topological
entropy. The procedure also gives some detailed information on the
support of the measure of maximal entropy.}}
\vskip .4truein
\centerline{\vbox{\hsize 4.5in \noindent Keywords: golden mean subshift,
hard square gas, lattice action, topological entropy, measure of maximal entropy \hfill\break
AMS Classification: 52C15, 58F08, 60K35, 82C20 \hfill\break
Running head: Hard Core entropy}}

\vfill
\eject

\hfill {}\break
\noindent {\subtitle Introduction}
\vskip .3truein
\noindent The {\bf Hard Core Model/Golden Mean Subshift/Independent Sets} is a highly useful model in various disciplines as witnessed by its many appearances under distinct names in fields like Statistical Mechanics/Symbolic Dynamics/Theoretical Computer Science respectively. Perhaps its most intuitive description is that of the interaction of identical hard particles whose only interaction is repulsion on contact. It has been studied for decades in various discrete set ups with some notable break throughs like Baxter's solution on the triangular lattice ([B1]) yet its general treatment still seems elusive.

The model can be considered in different regimes. In the {\bf high density/low temperature case} the characterization of the allowed configurations is essentially a packing problem ([E]). Here we concentrate on the {\bf loose packing/high temperature/entropic regime} where configurations don't have such long range order and hence no striking geometric structure. The key unanswered question is the exponential size of the configuration space which in turn boils down to computation of the topological entropy. For almost all lattice set-ups the exact answer in not known. Here we try to alleviate the situation a little bit by establishing a procedure that in principle yields arbitrarily good estimates of the entropy on the lattices and more over also gives some insight to the structure of the measure of maximal entropy.

\vskip .5truein
\noindent {\subtitle 1. Set-up}
\vskip .3truein
\noindent Let ${\rm {\bf L}}$ be a planar lattice. Subsequently we will consider mostly one of the regular lattices (square (${\rm {\bf Z}}^2$), honeycomb (${\rm {\bf H}}$) or triangular (${\rm {\bf T}}$)). In a few cases we illustrate the principles developed on more exotic stages like square lattice with the Moore neighborhood (${\rm {\bf Z}}^2{\rm M}$, nearest neighbors are within hop count one and two) or the Kagom\'e lattice (${\rm {\bf K}}$).

A configuration on ${\rm {\bf L}}$ satisfying the {\bf Hard Core Rule} is an element in $X=\{0,1\}^{{\rm {\bf L}}}$ where no two $1$'s can be nearest neighbors. This rule can naturally be viewed as a zero range infinite repulsive potential i.e. a hard exclusion rule not unlike that in hard sphere packing. Call the collection of configurations $X_{\rm {\bf L}}^{hc}.$

The exclusion rule naturally imposes a sublattice split on ${\rm {\bf L}}$ if it is a $k$-partite graph. For example on ${\rm {\bf Z}}^2$, a bipartite graph, one can man all sites on ${\rm {\bf 2Z}}^2$ (${\rm {\bf Z}}^2$ rescaled by $\sqrt{2}$ and rotated by $45^\circ$) with $1$'s and the rest of ${\rm {\bf Z}}^2$ must then be all $0$'s. Call the former the  {\bf even sublattice}, ${\rm {\bf L}}_e$ and the latter the {\bf odd sublattice}, ${\rm {\bf L}}_o$ (it is a $(1/2,1/2)$-shifted copy of the former). In rendering these we will present the even/odd sublattices as {\bf circle/dot sublattices}. In a similar fashion ${\rm {\bf H}}$ splits into two identical sublattices and ${\rm {\bf T}}$ (a tripartite graph) into three \lq\lq thinned\rq\rq\ copies of  ${\rm {\bf T}}.$ Both in the dense packing regime of [E] and in the the loose packing, entropic regime of this paper, this splitting will be highly relevant.

\vskip .3truein
\noindent Let $X_0\subset X.$ The standard measure of variation in the assignments on the lattice sites is as follows (as usual $|\cdot|$ means cardinality and $x|_A$ means the restriction of $x$ onto set $A$):

\proclaim Definition: The {\bf topological entropy} of the set $X_0$ is
$$h^{X_0}_{top}=\lim_{n\rightarrow\infty}{1\over n}\ln{\left|\left\{x|_{A_n}\ |\ x\in X_0\right\}\right|}$$
where $|A_n|=n$ and the sequence $\{A_n\}$ grows in a sufficiently regular fashion.\par

\noindent {\bf Remark:} For instance for the full shift on any lattice $h^X_{top}=\ln{2}$ (indicating two independent choices per lattice site). If $L={\rm {\bf Z}}$ the hard core model is explicitly solvable and a standard transfer matrix argument implies that $h^{\rm {\bf Z}}_{top}=\ln{\left({{1+\sqrt{5}}\over 2}\right)}$ (see e.g. [W]). For two and higher dimensional lattices the matrix argument breaks down and the exact value of the hard core topological entropy remains an unsolved problem. In this paper we try to approach and in particular approximate it in a novel way.

\vskip .2truein
\noindent From the general theory of lattice dynamical systems (e.g. [W]) it is known that shift invariant probability measures on a space of configurations, $\cal M$, satisfy the maximum principle
$$h_{top}=\sup_{\cal M}h_\mu$$
where $h_\mu$ is the measure-entropy. The special measures yielding the equality are {\bf measures of maximal entropy}. For two and higher dimensional systems they are in general not unique. In all our subsequent cases they are believed to be so, but we do not actually need this knowledge.

\vskip .2truein
\noindent Note that if all the limits exist, one could write $$h_{top}={1\over p}\sum_{i=1}^p h_{top}^{(i)}\leqno{(1.1)}$$ where $h_{top}^{(i)}$ is the topological entropy of the configurations on the $i^{\rm th}$ sublattice when the original lattice partitions into $p$ sublattices. However these entropies cannot in general be computed independently but rather depend on each other heavily like in the case of hard core. However it is still possible to imitate (1.1) by introducing a sequential approximation of the measure of maximal entropy $\mu.$

\vskip .5truein
\noindent {\subtitle 2. Lower bounds}
\vskip .3truein

\noindent We now proceed to establish lower bounds for the topological entropy using the sublattice partition representation (1.1) and a sequential fill-in scheme to overcome the dependencies. To keep the ideas clear we first treat the case of the hard core rule splitting the lattice to two sublattices and only after that generalize.

\vskip .2truein
\noindent Let $N_e$ denote an all-0 nearest neighbor neighborhood of a site on the odd lattice in the even lattice. In the case of ${\rm {\bf Z}}^2$ lattice the sites in $N_e$ form the vertices of an {\bf even unit diamond}, $\diamondsuit_e$ ($\diamondsuit_o$). On the honeycomb and triangular lattices these sites form triangular arrangements, $\bigtriangleup$ or $\bigtriangledown$ or a hexagon.

It will become quite useful to think the fill-in in terms of forming a tiling. The pieces are 0/1-{\bf tiles} which in ${\rm {\bf Z}}^2$ case are either 0/1-diamonds (as above) depending on whether the center site carries 0 or 1. On the hexagonal and triangular lattices the tiles are 0/1-{\bf (unit) hexes}. Once a sublattice is chosen, one can tile the plane using any combination of 0/1-tiles centered on the sublattice vertices.
 
Recall that the Bernoulli measure with parameter $p$, $B(p)$, assigns 1's independently with probability $p$ to each (sub)lattice site and 0's otherwise. Its entropy, denoted by $h_{B}(p),$ is $-p\ln{p}-(1-p)\ln{(1-p)}.$

\proclaim Proposition 2.1.: The topological entropy of the hard core model on a lattice with a two-way sublattice split is given by
$$h_{top}={1\over 2}\Big\{h_{top}^{(e)}+{\rm {\bf P}}\left(N_e\right)\ln{2}\Big\}\ ,\leqno{(2.1)}$$
where $h_{top}^{(e)}$ is the entropy of the measure of maximal entropy computed from to the even sublattice alone.\par

\vskip .2truein
\noindent {\bf Proof:} The representation (2.1) follows from (1.1) by observing that the maximum entropy is obtained by first assigning the marginal of the measure of maximal entropy to the even lattice and then filling in the non-blocked sites on the odd lattice. These are centered at the even unit diamonds. The non-blocked sites must be filled with $B(1/2)$ to obtain the maximal entropy on the odd lattice, hence the factor $h_{B}(1/2)=\ln{2}$. \hfill\QED

\vskip .2truein
\noindent The principle in the Proposition can be directly applied to square and honeycomb lattices. A further argument is required to cover all regular lattices. In the following result we present these arguments and further extend to Kagom\'e lattice, ${\rm {\bf K}}$ (a tripartite graph), as well as to the square lattice with Moore neighborhood, ${\rm {\bf Z}}^2{\rm M}$ (eight nearest neighbors, a 4-partite graph).

\proclaim Theorem 2.2.: The topological entropy of the hard core model is bounded from below on the square ($m=4$) and honeycomb ($m=3$) lattices by
$${\underline h}_{{\rm {\bf Z}}^2/{\rm {\bf H}}}(p)={1\over 2}\Big\{h_{B}(p)+(1-p)^m\ln{2}\Big\}\ ,\leqno{(2.2)}$$
on triangular ($m'=3$) and Kagom\'e ($m'=2$) lattices by
$${\underline h}_{{\rm {\bf T}}/{\rm {\bf K}}}(p,q)={1\over 3}\Big\{h_{B}(p)+(1-p)^{m'}[\ h_{B}(q)+[1-(1-p)q]^2\ln{2}\ ]\Big\} \leqno{(2.3)}$$
and on ${\rm {\bf Z}}^2{\rm M}$ lattice by
$$\eqalign{{\underline h}_{{\rm {\bf Z}}^2{\rm M}}(p,q,r)={1\over 4}\Big\{h_{B}(p)+&(1-p)^2\Big[\ h_{B}(q)+[1-(1-p)q]^4\ h_{B}(r) \cr 
+&(1-p)^2(1-q)^2\big[1-(1-(1-p)q)^2 r\big]^2\ln{2}\ \Big]\Big\}\ ,\cr}\leqno{(2.4)}$$
where $p,\ q$ and $r \in (0,1).$ \par

\vskip .2truein
\noindent {\bf Proof:} The lower bounds (2.2) follow simply from (2.1) of Proposition 2.1. by assigning $B(p)$ to the even sublattice since then ${\rm {\bf P}}\left(N_e\right)=(1-p)^{|N_e|}$ where the exponent is the number of elements in $N_e$ in ${\rm {\bf Z}}^2$ and $\rm {\bf H}$ respectively.

On the triangular lattice the sublattice split is three way. We call the parts the dot, circle and triangle sublattices. They are filled in three stages in the order $\circ\rightarrow\bullet\rightarrow\triangleright$. See Figure 1a and b for the notation and arrangement of the sublattices in a neighborhood of a triangle site.

Suppose the three sublattices are initially all empty. First fill-in the circle lattice with $B(p)$, hence the entropy contribution ${1\over 3}h_{B(p)}.$ Then fill-in all dot sites centered at $\bigtriangledown$ with $B(q)$, this implies the entropy increase ${1\over 3}(1-p)^3h_{B(p)}$ from the dot lattice.

To update the center site which is a triangle we need to know that its value is not forced. Hence
$$\eqalign{ &{\rm {\bf P}}({\rm center\ triangle\ not\ forced\ by\ nearest\ neighbor\ circle\ or\
 dot})\cr
&={\rm {\bf P}}({\rm no\ 1's\ in\ the\ hexagon\ of \ nearest\ neighbors\ of\ the\ triangle})\cr
&={\rm {\bf P}}({\bigtriangleup}={\underline 0}\ {\rm and}\ {\bigtriangledown}={\underline 0})={\rm {\bf P}}(c_2=c_4=c_5=0\ {\rm and}\ d_1=d_2=d_3=0) \cr
&={\rm {\bf P}}(d_1=d_2=d_3=0\ |\ c_2=c_4=c_5=0)\ {\rm {\bf P}}(c_2=c_4=c_5=0) \cr
&={\rm {\bf P}}(d_1=d_2=d_3=0\ |\ c_2=c_4=c_5=0)\ (1-p)^3 \cr
&=\left[\ {\rm {\bf P}}(d_1=0\ |\ c_2=c_4=c_5=0)\ \right]^3\ (1-p)^3 \cr
&=\left[\ {\rm {\bf P}}\left(c_1=1\ {\rm or}\ \{c_1=0\ {\rm and}\ d_1=0\}\ |\ c_2=c_4=c_5=0\right)\ \right]^3\ (1-p)^3 \cr
&=[p+(1-p)(1-q)]^3\ (1-p)^3\cr}\leqno{(2.5)} $$ 
which together with the choice $B(1/2)$ on the non-blocked dots gives (2.3).

The Kagom\'e lattice argument is similar to the triangular one. There are three sublattices involved, all identical copies of the Kagom\'e, only thinned and reoriented. For the nearest neighbors of a triangle-site see Figure 1c. Again we fill in the order $\circ\rightarrow\bullet\rightarrow\triangleright.$ In the last stage the probability of the triangle site being unforced is now

$$\eqalign{ &{\rm {\bf P}}(c_2=c_4=0\ {\rm and}\ d_1=d_2=0) \cr
&=\left[\ {\rm {\bf P}}\left(c_1=1\ {\rm or}\ \{d_1=0\ {\rm and}\ c_1=0\}\ |\ c_2=c_4=0\right)\ \right]^2\ (1-p)^2 \cr
&=[p+(1-p)(1-q)]^2\ (1-p)^2 \cr} $$

In the case of the square lattice with Moore neighborhood there is a four-way sublattice split. We denote and fill them in the following way: $\circ\rightarrow\bullet\rightarrow\triangleright\rightarrow\diamond$ (see Fig. 1d).

The two first terms of the formula (2.4) are straightforward since circles are laid independently and each dot has exactly two circle neighbors. Furthermore as above we can show that ${\rm {\bf P}}(\triangleright\ {\rm unforced})=(1-p)^2\left[p+(1-p)(1-q)\right]^4.$

For the diamond site at the center of Fig. 1d to contribute to the entropy we need to know the probability that it is unforced i.e. all entries in the punctured square $S$ rendered with dotted line in Fig 1d. are 0's:

$$\eqalign{ {\rm {\bf P}}\big(&S={\underline 0}\big) \cr
=&{\rm {\bf P}}\big({\rm all}\ \triangleright, \bullet\in S\ {\rm are\ 0's}\ |\ {\rm all}\ \circ\in S\ {\rm are\ 0's}\big)\ (1-p)^4 \cr
=&{\rm {\bf P}}\big(\triangleright\in S\ {\rm are\ 0's}\ |\ \circ, \bullet\in S\ {\rm are\ 0's}\big)\ (1-p)^4(1-q)^2 \cr
=&{\rm {\bf P}}\big(d_1=1\ {\rm or}\ \{d_1=d_2=0\ {\rm and}\ t_1=0\}\ |\ \circ, \bullet \in S\ {\rm are\ 0's}\big)\ (1-p)^4 (1-q)^2 \cr
=&\big[\ {\rm {\bf P}}\big(d_1=1\ |\ \circ, \bullet\in S\ {\rm are\ 0's}\big) \cr
&+{\rm {\bf P}}\big(d_1=d_2=0\ {\rm and}\ t_1=0)\ |\ \circ, \bullet\in S\ {\rm are\ 0's}\big)\big]^2\ (1-p)^4 (1-q)^2 \cr }\leqno{(2.6)}$$

One can compute the two probabilities in the last expression to be
$$2p(1-p)q+(1-p)^2(2-q)q$$
and
$$\left[p^2+2p(1-p)(1-q)+(1-p)^2(1-q)^2\right](1-r)$$
respectively. From these the formula in the square brackets in (2.6) can finally be simplified to the form $1-[1-(1-p)q]^2r.$ \hfill\QED

\vskip .4truein
\centerline{\hbox{
 \psfig{figure=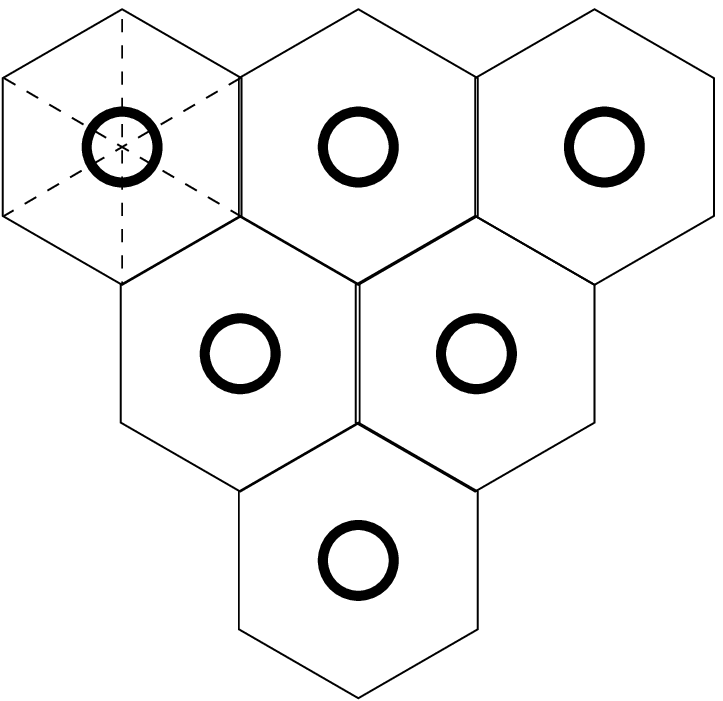,height=.9in}
 \hskip .3truein
 \psfig{figure=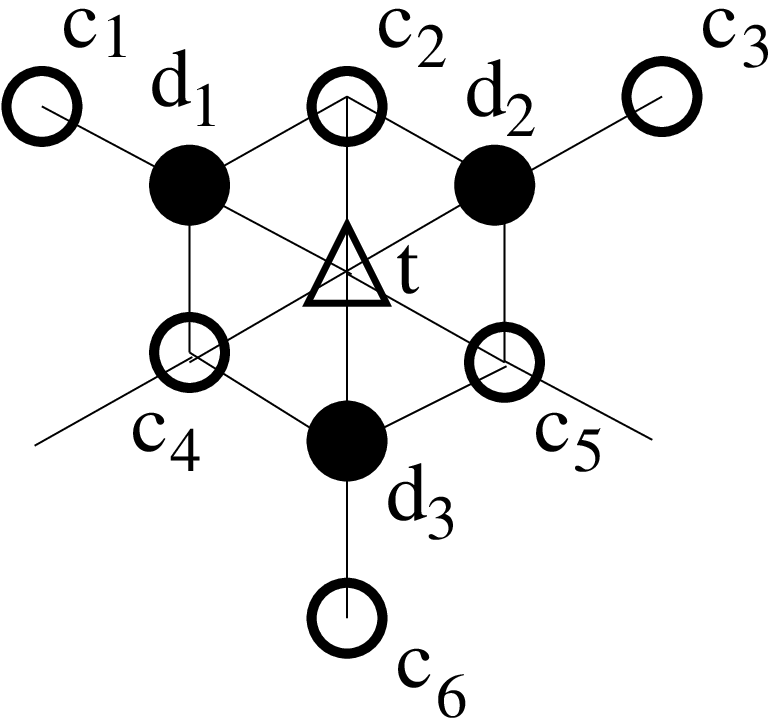,height=.9in}
\hskip .3truein
 \psfig{figure=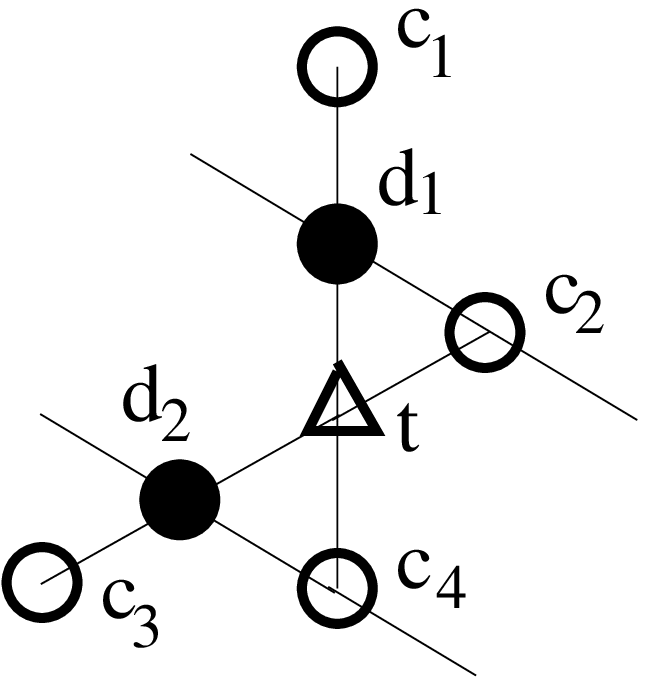,height=.9in}
 \hskip .3truein
 \psfig{figure=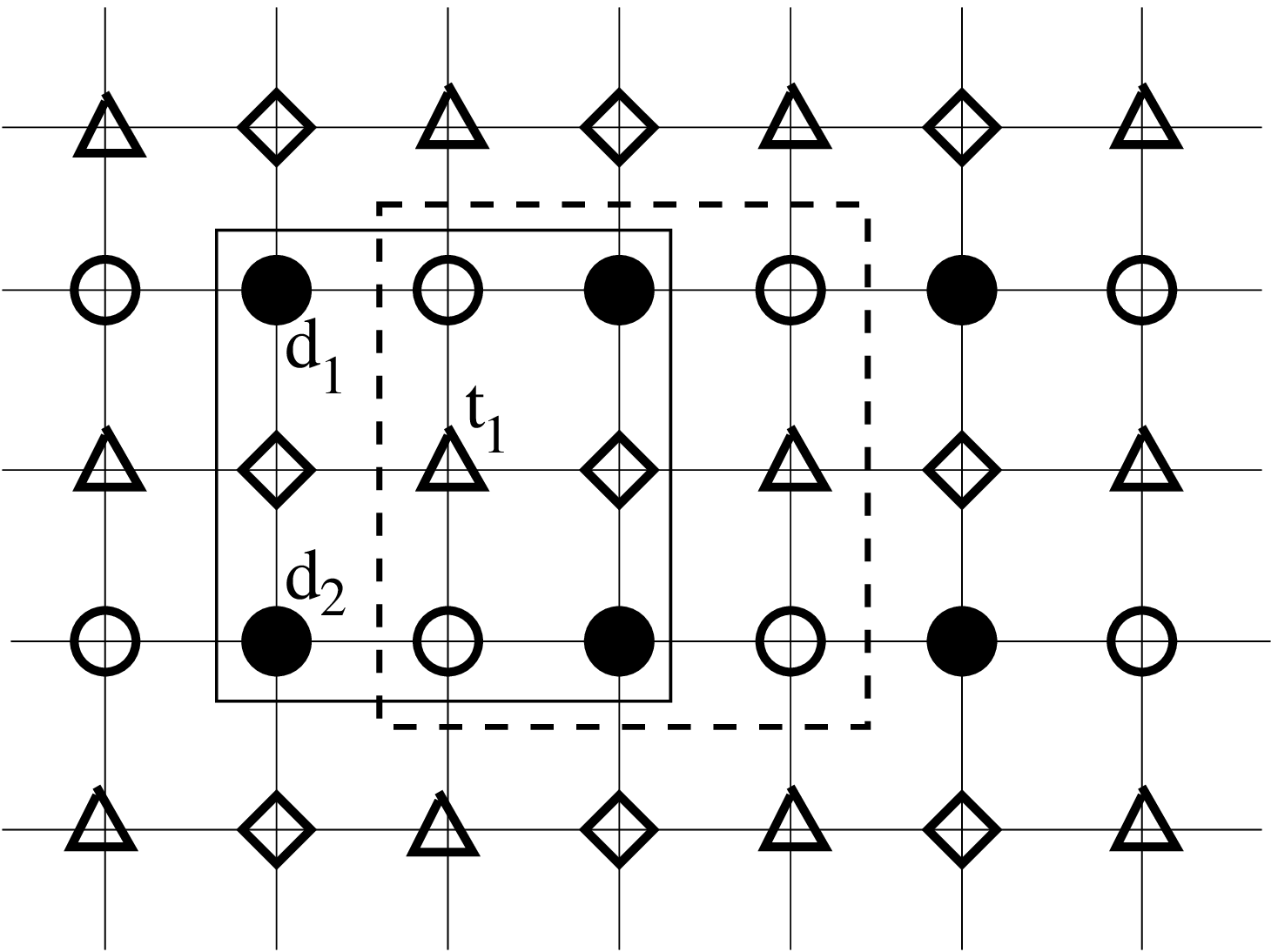,height=.9in}
}}
\vskip .2truein
\nobreak
\noindent {\bf Figures 1a, b, c and d.} Tiling of the triangular lattice with hexagonal 1-tiles and neighborhoods in triangular, Kagom\'e and ${\rm {\bf Z}}^2{\rm M}$ cases.
\vskip .3truein

\noindent The entropy bounds in (2.2) - (2.4) can be maximized with respect to the parameters using standard optimization routines in a desktop machine. 
 
\vskip .2truein

{
\offinterlineskip
\tabskip=0pt
\halign{ 
\vrule height2.75ex depth1.25ex width 0.6pt #\tabskip=1em &
#\hfil &\vrule \hfil # \hfil &  #\hfil &\vrule # &\hfil # \hfil &\vrule # & \hfil #\hfil & #\vrule width 0.6pt \tabskip=0pt\cr
\noalign{\hrule height 1pt}
& ${\rm {\bf L}}$ && max ${\underline h}_{\rm {\bf L}}$ && sublattice densities  && best estimates &\cr
\noalign{\hrule height 1pt}
& ${\rm {\bf Z}}^2$ && 0.3924 && (0.1702, 0.2370) && $0.4075\ (0.2266)_{\rm [MSS],[B2]}$ &\cr
\noalign{\hrule}
& ${\rm {\bf H}}$ && 0.4279 && (0.2202, 0.2371) && $0.4360\ (0.2424)_{\rm [B2]}$ &\cr
\noalign{\hrule}
& ${\rm {\bf T}}$ && 0.3253 && (0.1457, 0.1559, 0.1517) && $0.3332\ (0.1624)_{\rm [B2]}$ &\cr
\noalign{\hrule}
& ${\rm {\bf K}}$&& 0.3826 && (0.1944, 0.1948, 0.1866) &&  &\cr
\noalign{\hrule}
& ${\rm {\bf Z}}^2{\rm M}$ && 0.2858 && (0.119, 0.127, 0.130, 0.126) &&  &\cr
\noalign{\hrule height 1pt}
}}

\vskip .3truein
\noindent {\bf Table 1.} First lower bounds for Hard Core topological entropy and the corresponding sublattice densities for some 2-d lattices. To the right we have indicated the best numerical estimates for the entropy and corresponding density (in parenthesis) found in the literature.

\vskip .2truein
\noindent We note that while the topological entropy has been computed in the square lattice case to a great accuracy (e.g. in [B2] to some 40 decimal places) the corner transfer matrix methods used in these numerical studies attack the problem in a very different way. Our aim is not to compete in decimal count but rather present an alternative method applicable in many lattice set-ups to estimate the entropy which simultaneously yields some explicit information on the generic configurations/the measure of maximal entropy.

\vskip .2truein
\noindent The measure of maximal entropy doesn't need to be unique for a 2-d lattice model but in the case of hard square gas it is. This follows from the Dobrushin criterion ([DS], [RS]). Using this knowledge and the results above we now establish bounds for the density of 1's in the generic configurations. The exact value of the upper bound in the following result is in the Proof but we prefer to give the statement in this more explicit form.

\proclaim Proposition 2.3.: In the square lattice case the density of 1's at the equilibrium is in the interval $(0.21367,0.25806)$. \par

\noindent {\bf Proof:} Let $\rho_e$ be the density of 1's on the even lattice and let $c$ denote the expected number of 0's that a 1 forces on the odd lattice. Since exactly half of the non-forced sites will be 1's it must by the uniqueness of the measure of maximal entropy hold that $(2+c)\rho_e=1.$ Hence under
it
$${\rm {\bf P}}(x_i=0)={{1+c}\over{2+c}}\qquad 
{\rm {\bf P}}(x_i=1)={1\over{2+c}}\qquad
{\rm {\bf P}}\left(\diamondsuit_e\right)={2\over{2+c}}$$
on both lattices. $\diamondsuit_e$ is the 0-diamond as defined in the beginning of the section.

The entropy of any distribution on the even lattice with 1-density $\rho_e$ is bounded from above by the entropy of the Bernoulli distribution with parameter $\rho_e.$ Hence the total entropy at that 1-density is bounded from above by
$${1\over 2}\left(h_{B}\left({1\over{2+c}}\right) + {2\over{2+c}}\ln{2}\right).$$ 
This expression bounded by $h_{top}$ of Table 1 ([B2] or [MSS]) yields an upper bound for $c$, $2.6801$ which in turn gives the lower bound for $\rho_e$.

The upper bound for $\rho_e$ follows from  a lower bound for $c$ which we establish using a monotonicity argument. The 1's on say the even lattice are $B(1/2)$-distributed on the non-forced sites. Call this set $F$ and pick a site on it which has a 1. How many sites will this entry block? Let $F'$ be a superset of $F.$ Then clearly ${\rm {\bf E}}(c|\ F)\ge {\rm {\bf E}}(c|\ F')$ as in a bigger domain the 1 is more likely to share the blocking with a nearest neighbor 1 on the same sublattice. Hence a lower bound is obtained by calculating the blocking for a 1 with its eight nearest neighbors also in $F$.  Enumerating the $2^8$ possible neighborhood configurations and weighting them uniformly according to the $B(1/2)$-distribution we get the lower bound for $c$: $15/8.$ This in turn implies the upper bound for $\rho_e$, $8/31$.
\hfill\QED

\vskip .2truein
\noindent Since our first estimate for the lower bound on ${\rm {\bf Z}}^2$ is associated with densities incompatible with Proposition 2.3. we will try out a symmetric variant of the theme. The (near) equality of the densities on the sublattices should be a natural property of a measure corresponding to a good lower bound since the measure of maximal density is believed to be unique in all our cases. In the last three cases in Table 1. the non-equality of the densities isn't far off but for the first two we present an \lq\lq equalization\rq\rq.
 
\proclaim Proposition 2.4.: To achieve equal densities of 1's on each of the sublattices one needs to replace the $B(1/2)$ distribution in the last stage of the measure construction by $B(p')$ and thereby $\ln{2}$ in (2.2) by
$h_{B(p')}$, where $p'=p(1-p)^{-|N_e|}.$ \par

\noindent {\bf Proof:} In the case of two sublattices after $B(p)$ distribution of 1's on the even lattice there are a density of $(1-p)^{|N_e|}$ unforcing neighborhoods on this sublattice. These have to produce the correct density of 1's on the odd lattice, hence we need the even lattice flip probability $p'$ to satisfy $p'(1-p)^{|N_e|}=p.$ \hfill\QED

\vskip .2truein
\noindent Using the Proposition one can optimize the square lattice topological entropy bound to (a slightly worse value) 0.3921 at joint density level 0.2015. In view of Proposition 2.3. this indicates that the entropy generating 1's are not yet packed in densely enough. In the case of the honeycomb lattice the corresponding values are 0.427875 at 0.2284.


\vskip .5truein
\noindent {\subtitle 3. Higher order blocks}
\vskip .3truein
\noindent

\noindent To improve the entropy bounds and more importantly to get some insight into the character of the measure of maximal entropy we now consider more complicated optimization schemes involving Bernoulli-distributed blocks on sublattices. We first illustrate the ideas on hexagonal and triangular lattices.

A {\bf three-hex} is a obtained by gluing together three unit hexes so that each has two joint sides. Figure 2a. illustrates three such three-hexes next to each other (for reference lattice edges are indicated as thin dotted lines in one of the unit hexes). Note that the unit tiles on each of them are all centered on the same sublattice, the circle lattice in this case (call the tile a {\bf circle three-hex}). The dots of the other sublattice are all in the centers of the three-hexes or on their extremities (three of them are indicated). Three-hexes of the same orientation obviously tile the plane.

\vskip .4truein
\centerline{\hbox{
 \psfig{figure=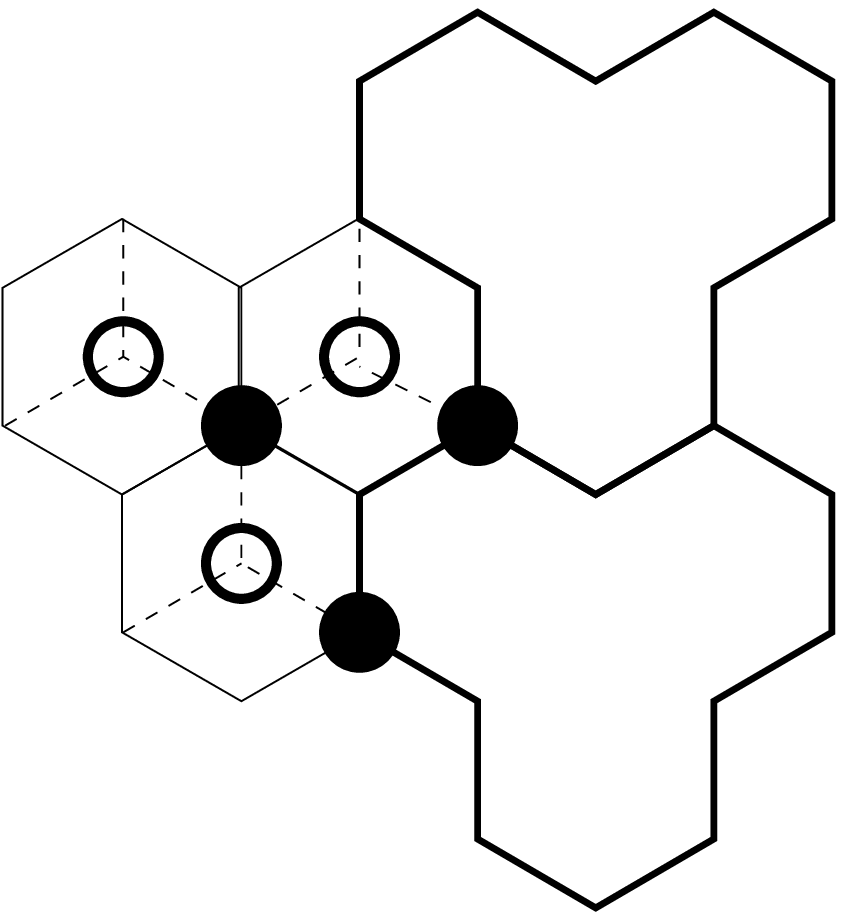,height=1.3in}
 \hskip .4truein
 \psfig{figure=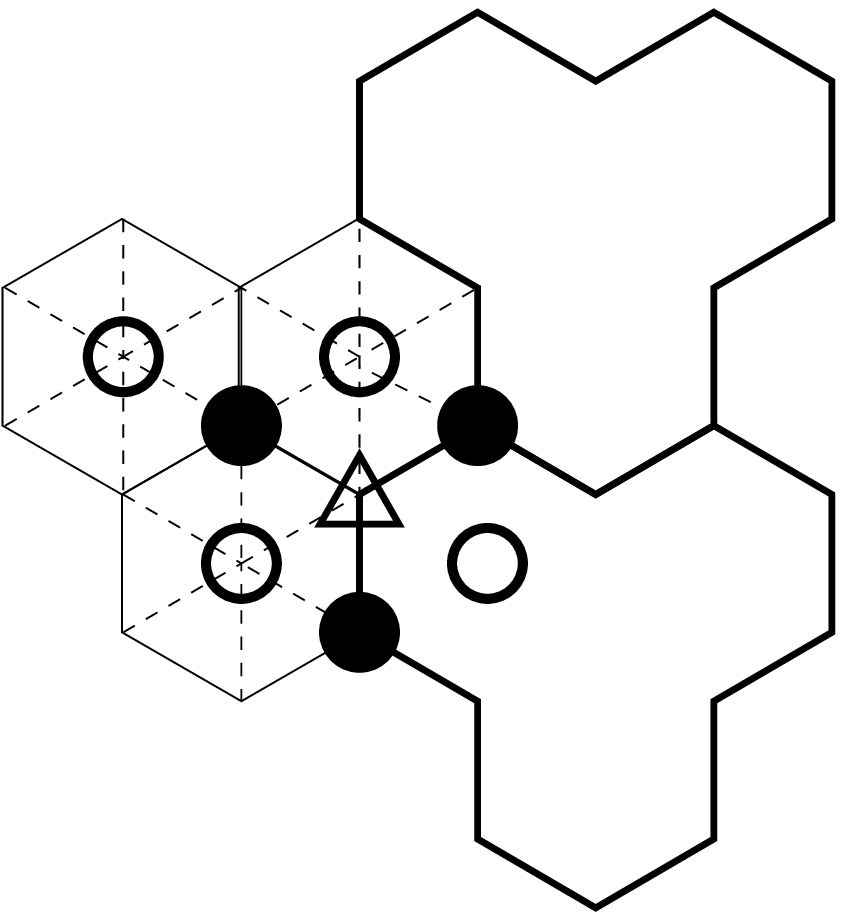,height=1.3in}
 \hskip .8truein
 \psfig{figure=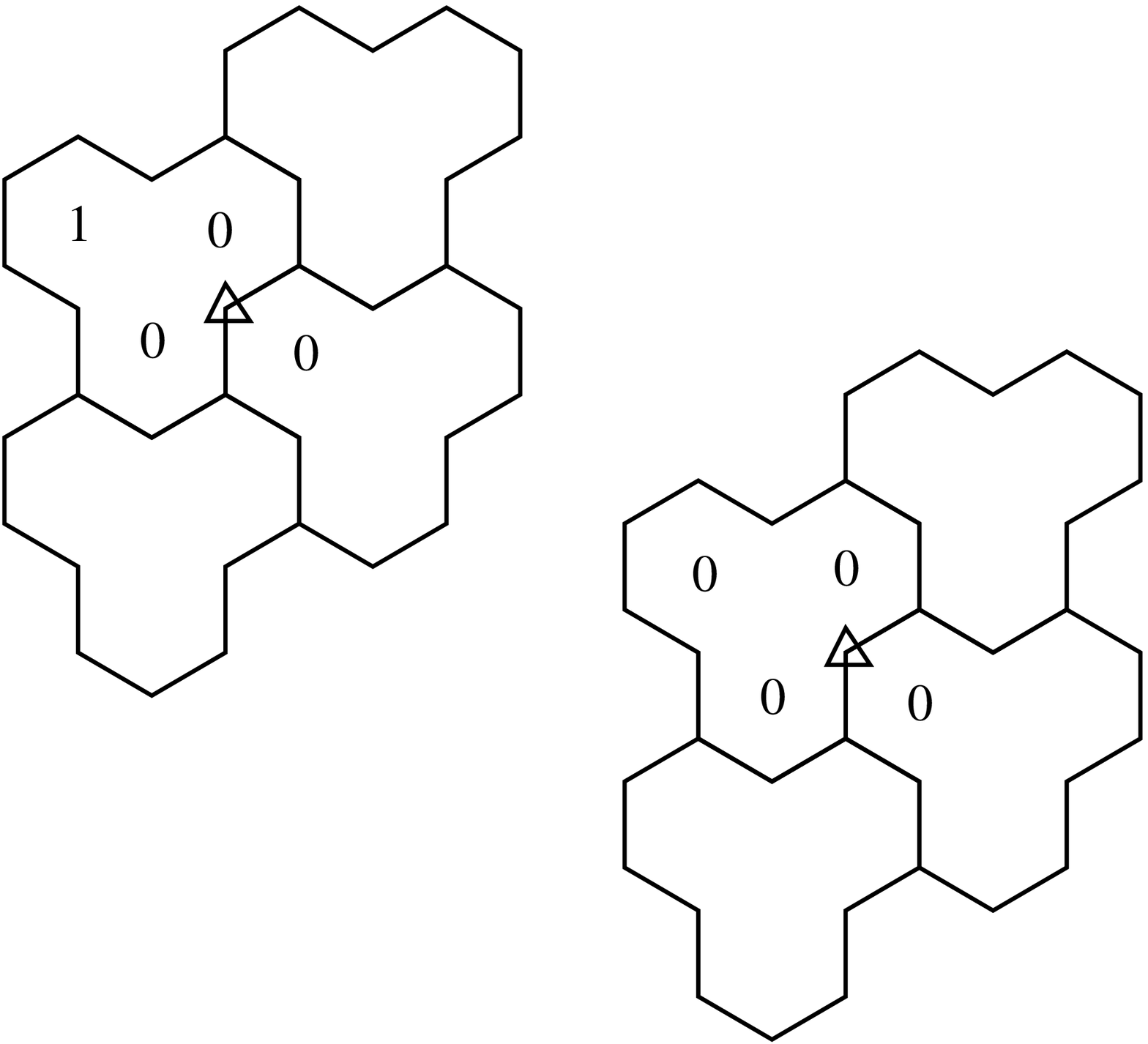,height=1.3in}
}}
\nobreak
\vskip .2truein
\noindent {\bf Figures 2a, b and c.} 3-hex arrangements in hexagonal and triangular cases.
\vskip .3truein

\noindent Let $B({\bf p}),\ {\bf p}=(p_0,p_1,p_2,p_3)$ be the Bernoulli distribution on circle three-hexes with the probability that the three-hex has exactly $k$ 1-tiles in it in a given orientation being $p_k$ (so $p_0+3p_1+3p_2+p_3=1).$ Its entropy is then $h_{B}^{(3)}({\bf p})=-p_0\ln{p_0}-3p_1\ln{p_1}-3p_2\ln{p_2}-p_3\ln{p_3}.$

\proclaim Theorem 3.1.: Let $a({\bf p})=p_0+2p_1+p_2.$ For the hexagonal lattice the Hard Core entropy is bounded from below by
$${\underline h}_{\rm {\bf H}}^{(3)}({\bf p})={1\over 6}\Big\{h_{B}^{(3)}({\bf p})+\left[p_0+2a({\bf p})^3 \right]\ln{2}\Big\} \leqno{(3.1)}$$
\noindent and for the triangular lattice a corresponding bound is
$$\eqalign{ {\underline h}_{\rm {\bf T}}^{(3)}({\bf p},q)={1\over 9}\Big\{h_{B}^{(3)}({\bf p})&+\left[p_0+2a({\bf p})^3 \right]h_B(q) \cr
&+3\left[p_1+p_0(1-q)\right]a({\bf p})^3(2-q)^2\ln{2}\Big\}\ ,\cr} \leqno{(3.2)} $$
where $p_i,\ q\in (0,1).$
\par

\vskip .2truein
\noindent {\bf Proof:} For the construction of the measure we will fill in the lattice in the order $\circ\rightarrow\bullet.$ If the circle three-hexes are distributed Bernoulli with parameter ${\bf p}$ the entropy contribution from the circle lattice will be
${1\over 2}{1\over 3}h_B({\bf p})$ 
where the factors result from the sublattice density and the fact that we distribute triples. As in Proposition 2.1. in the next stage the maximal entropy choice for the unforced sites on the dot lattice is the $B(1/2)$ distribution. The total density of sites available is computed at two different types of dot sites (as in Fig. 2a, the three dots indicated) and is
${1\over 3}\left[p_0+2\left(p_0+2p_1+p_2\right)^3\right]$
where the coefficient 2 and the power 3 follow from the fact that at two of the three dot sites three adjacent three-hexes coincide. These formulas combined and simplified yield $(3.1)$.

On the triangular lattice a third sublattice enters and the fill-in order is then $\circ\rightarrow\bullet\rightarrow\triangleright.$ The entropy contribution from the Bernoulli circle three-hexes is now ${1\over 3}{1\over 3}h_B({\bf p})$ since each sublattice is identical, hence of density $1/3.$

In the second stage the unforced dot sites are filled with $B(q)$ distribution. Their density is computed as above to be
${1\over 3}\left[p_0+2\left(p_0+2p_1+p_2\right)^3\right],$ hence the entropy contribution from dot lattice will be this expression multiplied by ${1\over 3}B(q).$

In the final stage the unforced triangle sites are filled by $B(1/2).$ Their density in the full lattice is
$$\eqalign{ {1\over 3}{\rm {\bf P}}(&{\rm nearest\ neighbor\ } \circ {\rm \ and\ } \bullet {\rm \ sites\ all\ 0's}) \cr
=&{1\over 3}\Big\{p_1(p_0+2p_1+p_2)\left[(p_1+2p_2+p_3)+(p_0+2p_1+p_2)(1-q)\right]^2 \cr
&+p_0(p_0+2p_1+p_2)(1-q )\left[(p_1+2p_2+p_3)+(p_0+2p_1+p_2)(1-q)\right]^2\Big\}, \cr }
\leqno{(3.3)}$$
which results from considering the two different arrangements of four neighboring three-hexes as shown in Fig. 2c. (top and bottom cases for the top and bottom expressions in $(3.3)$). The formulas merged and simplified result in $(3.2).$ \hfill\QED

\vskip .2truein

{
\offinterlineskip
\tabskip=0pt
\halign{ 
\vrule height2.75ex depth1.25ex width 0.6pt #\tabskip=1em &
#\hfil &\vrule \hfil # \hfil &  #\hfil &\vrule # &\hfil # \hfil &\vrule # & \hfil #\hfil & #\vrule width 0.6pt \tabskip=0pt\cr
\noalign{\hrule height 1pt}
& ${\rm {\bf L}}$ && max ${\underline h}_{\rm {\bf L}}$ && $(p_0,p_1,p_2,p_3)$,  $ q$  && sublattice densities &\cr
\noalign{\hrule height 1pt}
& ${\rm {\bf H}}$ && 0.4304 && (0.504, 0.110, 0.048, 0.021) && (0.2276, 0.2376)  &\cr
\noalign{\hrule}
& ${\rm {\bf T}}$ && 0.3265 && (0.64, 0.092, 0,025, 0.010), 0.25  && (0.153, 0.155, 0.151) &\cr
\noalign{\hrule height 1pt}
}}

\nobreak
\vskip .2truein
\noindent {\bf Table 2.} Optimized lower bounds and densities for three-hex Bernoulli blocks. 

\vskip .2truein
\noindent {\bf Remarks: 1.}  The Kagom\'e lattice case can be done in an identical fashion.\hfill\break
\noindent {\bf 2.} Note that apart from improvements in the entropy bounds, almost all of the sublattice densities have increased (in comparison to values in Table 1) indicating a better packing of the 1's on the sublattices. Moreover they have significantly less variation which is to be expected since the densities are equal for the measure of maximal entropy.

\vskip .3truein
\noindent Let us now return to our original motivation, the Hard Core on the square lattice. Compounding the principles above and some further ideas we will implement an increasing sequence of lower bounds converging to the topological entropy. Along the way we'll get more explicit information on the configurations favored by the measure of maximal entropy.

1-tiles in the ${\rm {\bf Z}}^2$ case are diamonds of side length $\sqrt{2}$ centered on either of the two sublattices. $k${\bf -omino} is formed by gluing together $k$ such 1-tiles along edges. If $k=n^2$ and the 1-tiles are in a diamond formation we call them a $n\times n$ {\bf -blocks}. There are $2^{n^2}$ of them. The optimization results in Section 2 were for the $1\times 1$-blocks.

Consider next $2\times 2$ -blocks. There are $16$ of them, but after assuming isotropy for them i.e. that blocks that are rotations of each other are distributed with equal probability (inevitable when measure of maximal entropy is unique), there are only five free parameters for Bernoulli distribution $B({\rm {\bf p}})$ on them (${\rm {\bf p}}=(p_0,p_1,p_{21},p_{22},p_3,p_4)$, $p_0+4p_1+4p_{21}+2p_{22}+4p_3+p_4=1.$ Here the first subindex of $p$ refers to the number of 1's in the block and $p_{22}$ and $p_{21}$ denote the two different arrangement of two 1's in the block (side by side and across)).

The entropy contribution from the even lattice (on which we distribute first the 1's using $B({\rm {\bf p}})$ is now
$$-{1\over 4}\Big\{p_0\ln{p_0}+4p_1\ln{p_1}+4p_{21}\ln{p_{21}}+2p_{22}\ln{p_{22}}+4p_3\ln{p_3}+p_4\ln{p_4}\Big\}\ .\leqno{(3.4)}$$
The density of the unforced sites on the odd lattice can be computed from the three cases indicated in Figure 3a. and results in 
$${1\over 4}\left\{p_0+2\left(p_0+2p_1+p_{21}\right)^2+\left(p_0+3p_1+2p_{21}+p_{22}+p_3\right)^4 \right\}\ .\leqno{(3.5)}$$
These combined yield a lower bound for $h_{top},$ which is optimized in Table 3 (second row).


\vskip .2truein
\centerline{\hbox{
 \psfig{figure=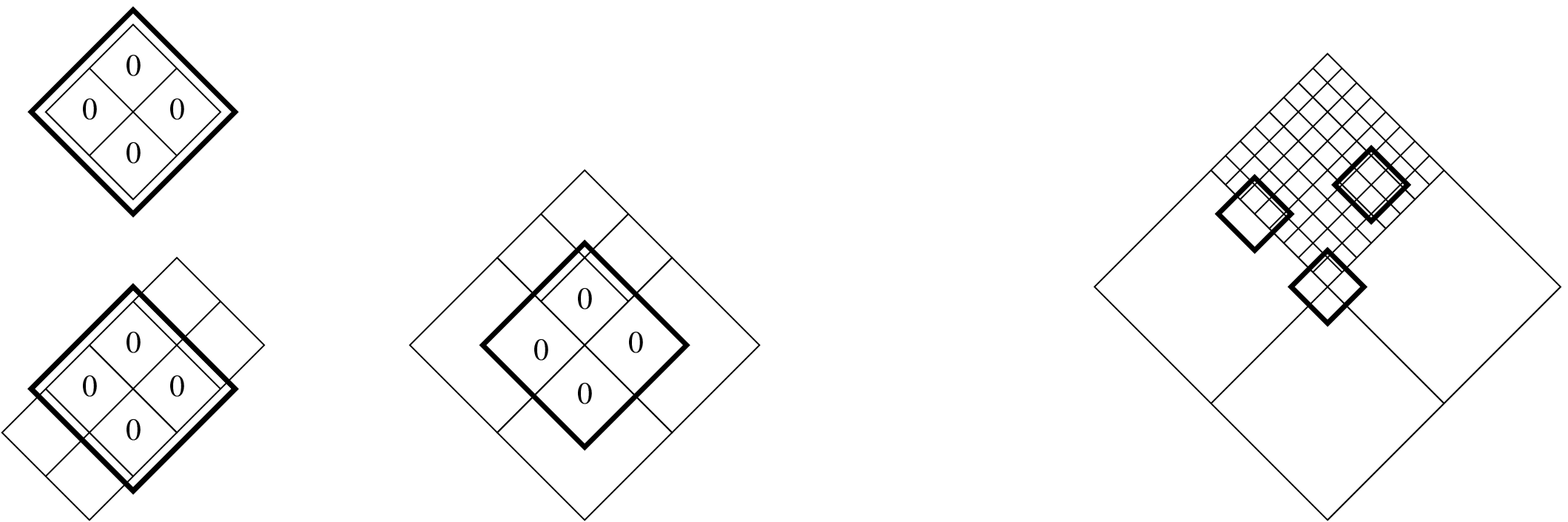,height=.9in}
 \hskip .7truein
\psfig{figure=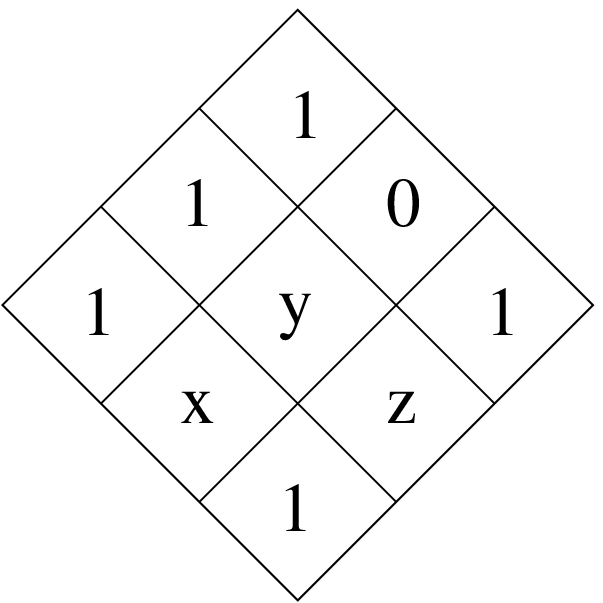,height=.6in}
 \hskip .3truein
 \psfig{figure=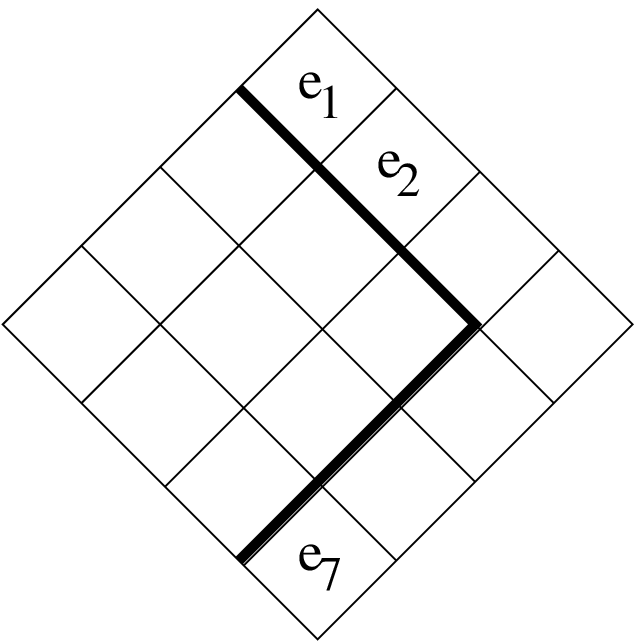,height=.6in}
}}
\vskip .1truein
\nobreak
\noindent {\bf Figures 3a {(\rm 3)}, b, c and d.} $n\times n$ -blocks and update window in ${\rm {\bf Z}}^2$ case. Reductions in a $3\times 3$ -block and the extension.  

\vskip .3truein
{
\offinterlineskip
\tabskip=0pt
\halign{ 
\vrule height2.75ex depth1.25ex width 0.6pt #\tabskip=1em &
#\hfil &\vrule \hfil # \hfil &  #\hfil &\vrule # &\hfil # \hfil &\vrule # & \hfil #\hfil & #\vrule width 0.6pt \tabskip=0pt\cr
\noalign{\hrule height 1pt}
& block size && max ${\underline h}_{{\rm {\bf Z}}^2}$ && sublattice densities  && red/init variables &\cr
\noalign{\hrule height 1pt}
& $1\times 1$ && 0.392421 && (0.1702, 0.2370) && 1  &\cr
\noalign{\hrule}
& $2\times 2$ && 0.39877 && (0.1993, 0.2254)  && 5/15 &\cr
\noalign{\hrule}
& $3\times 3$ && 0.4014 && (0.2073, 0.2254)  && 46/511 &\cr
\noalign{\hrule height 1pt}
}}

\vskip .2truein
\noindent {\bf Table 3.} Optimized lower bounds and densities for a few Bernoulli blocks on ${\rm {\bf Z}}^2.$

\vskip .2truein
\noindent In the block size $3\times 3$ there are initially 511 free block probabilities to optimize. When rotational invariance is imposed the variable number is reduced and additionally we will expect blocks that are reflections of each other to have equal probabilities at the optimum. After these two types of symmetries are accounted the number of free variables will be 101.


In this size and in larger blocks another feature appears which enables further variable weeding. Consider the block in Figure 3c. The symbol assignments in sites $x,\ y$ and $z$ are irrelevant in the sense that the existing 1's in the $3\times 3$ -block already force all the odd sites (to carry 0's) that $x,\ y$ and $z$ might force if any of the them were 1's. Hence there are $2^3$ blocks of equal probability. This combined with the symmetry assumptions above yield the total of 64 blocks with identical probabilities at the optimum (this is actually the maximum reduction achievable in this block size). Combing through the set of all blocks for this feature will result in reduction by a factor about 11 to the final set of 46 variables. Their optimal values have been computed and the results are in Table 3.

Subsequently we call sites like $x,y$ and $z$ above {\bf weak} with respect the rest of the given block. Only the corner sites of a block cannot ever be weak.

The procedure of variable reduction is highly useful since the above rotational and reflection symmetry search as well as the weak site identification can be automated. Moreover the reduction improves significantly at every stage: for example in the next block size of $4\times 4$ the initial variable number of 65.536 shrinks 66-fold to 991 final free variables.

Note also that the optima in block size $n\times n$ can be utilized as indicated in Figure 3d to initiate the search in the next larger block size. Once e.g. the $3\times 3$ subblock optimum probability is known, the added half frame $(e_1,\ldots,e_7)$ should be assigned $B(p)$ entries with $p$ computed from $3\times 3$ blocks. With tailored optimization routines one should be able to deal with several thousands of variables in the larger block sizes. All the optimizations here were done with non-specialized code using {\sl Mathematica}.

\vskip .3truein
\noindent The optimal block probabilities satisfy a useful monotonicity property, that we establish next. For this let $B_i,\ i=1,2$ be $n\times n$ -blocks, whose subsets of 1's we refer to as $B_i^{(1)}.$ There is a partial order on the blocks via $B_i^{(1)}$ using the ordinary set inclusion. Let the optimal probabilities for the blocks be ${\rm {\bf p}}=(p_0,p_1,p_2,p_3,\ldots,p_l),\ l=2^{n^2}$ (no reductions done yet and no particular order in the coordinates).

\vskip .1truein
\proclaim Theorem 3.2.: Given two blocks $B_1$ and $B_2$ with optimal lower bound probabilities $p_1$ and $p_2$, if $B_1^{(1)}\subset B_2^{(1)}$ then $p_1\ge p_2.$ If $B_2^{(1)}\setminus B_1^{(1)}$ contains only weak sites with respect to $B_1^{(1)}$ then $p_1=p_2$, otherwise $p_1>p_2.$ \par

\vskip .1truein
\noindent {\bf Proof:} The optimal lower bound is given by
${\underline h}({\rm {\bf p}})={1\over {n^2}}\left\{-\sum_i p_i\ln{p_i}+{\rm {\bf P}}(N_e)\ln{2}\right\}$
where $N_e$ is the even $2\times 2$ -diamond of all 0's as in Section 2. Let $B_i$ be such that $B_1^{(1)}\subset B_2^{(1)}$ and let $p_1=p+\epsilon$,\ $p_2=p-\epsilon$,\ $0\le|\epsilon|<p.$ Denote by $h_\epsilon({\rm {\bf p}})$ the lower bound with the given $p_1$ and $p_2.$ To prove the result we will consider the entropy variation under the probability change of the two blocks: $\Delta h_\epsilon({\rm {\bf p}})=h_\epsilon({\rm {\bf p}})-h_0({\rm {\bf p}}).$ More explicitly
$$\eqalign{\Delta h_\epsilon({\rm {\bf p}})={1\over {n^2}}\Big\{\big[&-(p+\epsilon)\ln{(p+\epsilon)}-(p-\epsilon)\ln{(p-\epsilon)}+2p\ln{p}\big] \cr
+\big[&{\rm {\bf P}}_{1,\epsilon}(N_e)-{\rm {\bf P}}_1(N_e) \cr
+&{\rm {\bf P}}_{2,\epsilon}(N_e)-{\rm {\bf P}}_2(N_e) \cr
+&{\rm {\bf P}}_{4,\epsilon}(N_e)-{\rm {\bf P}}_4(N_e)\big]\ln{2}\Big\}, \cr} \leqno{(3.6)}$$
where ${\rm {\bf P}}_{k,\epsilon}(N_e)$ and ${\rm {\bf P}}_k(N_e)$ are the $N_e$-diamond probabilities computed from the different arrangements involving $k=1,2$ or $4$ $n\times n$ -blocks as in Figures 3a and b, for the block probability choices $p\pm\epsilon$ or $p$ for both.

By $\ln{(1+x)}\approx x$ the first square bracket behaves for small $\epsilon$ like $c_1\epsilon^2,\ c_1<0.$ 

 If $B_2^{(1)}\setminus B_1^{(1)}$ contains only weak sites with respect to $B_1^{(1)}$ then the blocks $B_i$ allow exactly the same sites to flip on the odd lattice hence each of the three last lines in (3.6) vanishes. The sole contribution to $\Delta h_\epsilon({\rm {\bf p}})$ then comes from the first square bracket and since this is negative for small but nonzero $\epsilon,$ it must be that $p_1=p_2$ at the optimum.   

If $B_2^{(1)}\setminus B_1^{(1)}$ contains non-weak sites with respect to $B_1^{(1)}$ let us first assume that they force $k$ odd interior sites (recall that the odd sites are the vertices of the grids in Figure 4. There are $(n-1)^2$ such interior sites in a $n\times n$ -block). Let $m$ be the number non-forced odd interior sites over block $B_1.$ Then  
$$\eqalign{ {\rm {\bf P}}_{1,\epsilon}(N_e)-{\rm {\bf P}}_1(N_e)&=\Big(\ldots + {(p+\epsilon)m\over {(n-1)^2}}+{(p-\epsilon)(m-k)\over{(n-1)^2}}+\ldots\Big) \cr
&-\Big(\ldots + {pm\over{(n-1)^2}}+{p(m-k)\over{(n-1)^2}}+\ldots\Big)={k\epsilon\over{(n-1)^2}}, \cr}$$
where the dots refer to the contributions from the other blocks. All these terms cancel out, since the other block probabilities are identical.

If non-weak sites only force odd interior sites then by geometry of the set-up the two last lines in (3.6) are immediately zero. If $e$ extra odd edge, off-corner sites are forced, similar argument than above gives estimate
$(c_2+{e\epsilon\over{4(n-1)}})^2-c_2^2, c_2>0$ for ${\rm {\bf P}}_{2,\epsilon}(N_e)-{\rm {\bf P}}_2(N_e)$ so the next to last line in (3.6) has the first order behavior $c_3\epsilon,\ c_3>0.$ Some added bookkeeping yields ${\rm {\bf P}}_{4,\epsilon}(N_e)-{\rm {\bf P}}_4(N_e)=(c_4+l\epsilon/4)^4-c_4^4\approx c_5\epsilon,\ c_5>0$ ($l$ is the number of odd corners forced). 

The leading order estimates for the four terms in the square brackets in (3.6) together yield $c_1\epsilon^2+d\epsilon,\ c_1<0,\ d\ge 0.$ If there are non-weak sites in $B_2^{(1)}\setminus B_1^{(1)}$ with respect to $B_1^{(1)}$, then $d>0.$ Hence $p_1>p_2$ must prevail at the optimum. \hfill\QED

\vskip .3truein
\noindent {\bf Remarks: 1.} Intuitively the result says that if neither of two even blocks gives more subsequent choice on the odd lattice, for maximum entropy one should weight them equally. Otherwise one should favor the one giving more choice on the odd. \hfill\break
\noindent {\bf 2.} One can readily see some chains imposed by the order in Figure 4: $0\prec 12\prec 23\prec 31$ or $0\prec 11\prec 21/22\prec 33$  etc. The monotonicity can be utilized in limiting the number of $n\times n$ -blocks optimized for larger values of $n$ (dropping blocks with least probability as dictated by the Theorem and with least multiplicity (most symmetric)).

\vskip .4truein
\centerline{\hbox{
 \psfig{figure=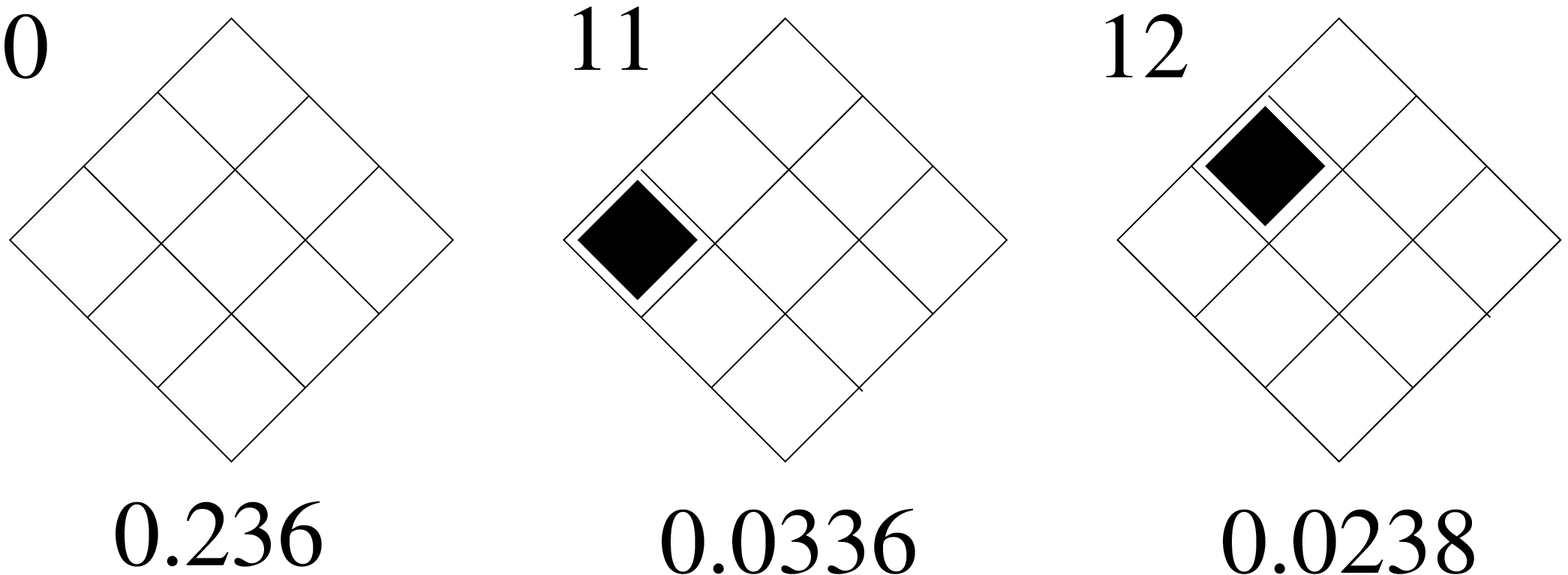,height=.55in}
\hskip .1truein
 \psfig{figure=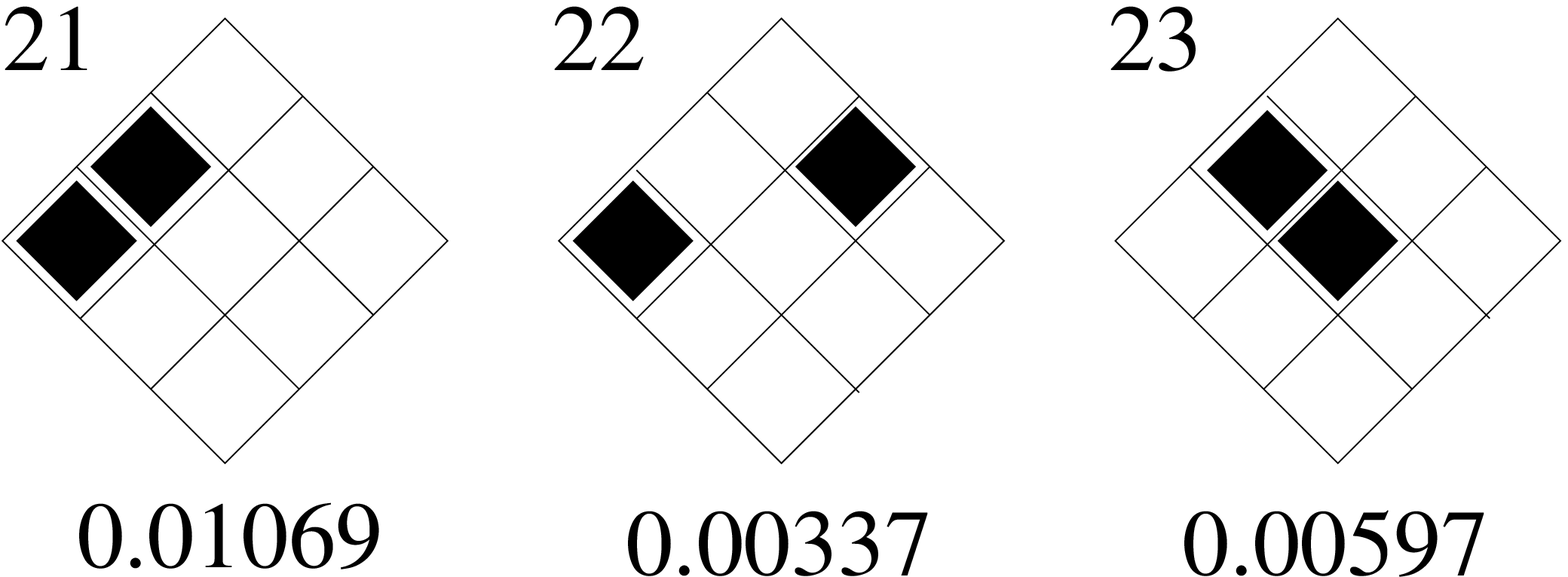,height=.55in}
\hskip .1truein
 \psfig{figure=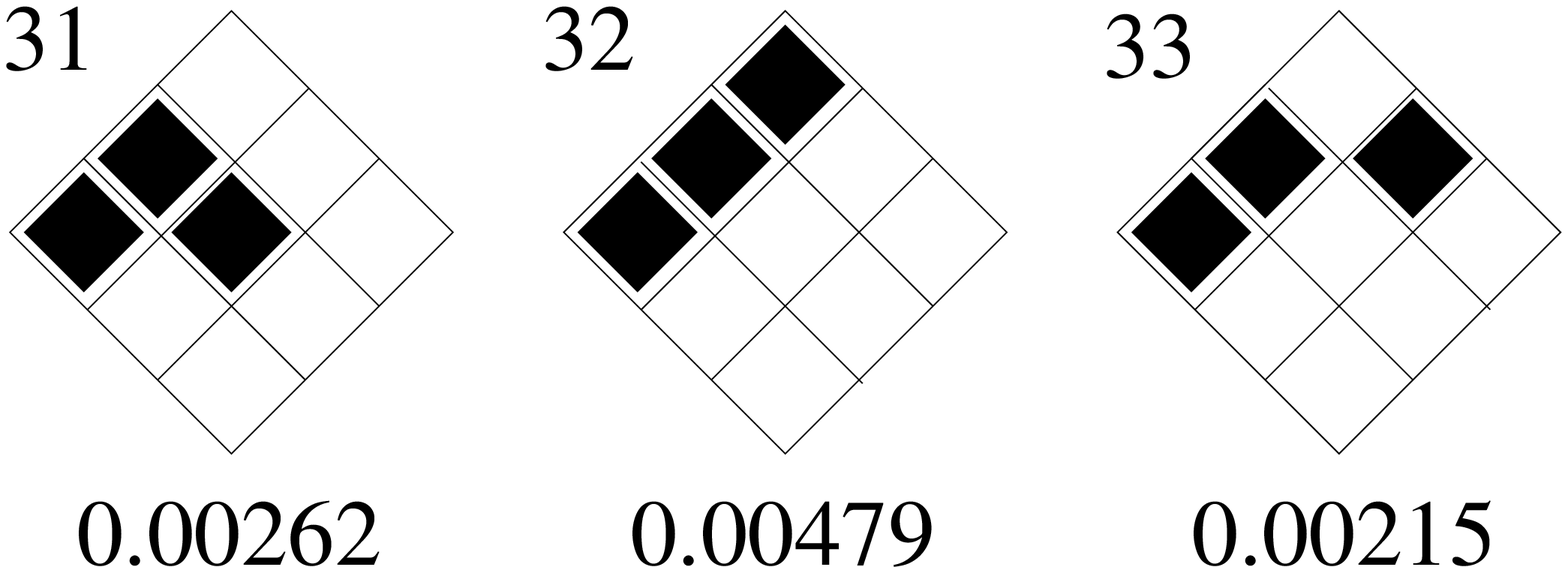,height=.55in}
}}
\vskip .1truein
\nobreak
\noindent {\bf Figure 4.} Prevalent $3\times 3$ -blocks with optimal probabilities without multiplicities.

\vskip .3truein
\noindent The correlation structure inside the measure of maximal entropy gradually presents itself in the Bernoulli approximations when we consider higher order blocks. Correlations between the blocks are zero because of independence, but within the blocks it is worth making comparisons.

By adding the optimum probabilities of all $3\times 3$ blocks at a given density level $k/9=0,1/9,\ldots ,1$ we obtain the \lq\lq density profile\rq\rq\ of this measure (here $k$ is the number of 1's in the block).   

Suppose next that we generate the $3\times 3$ blocks from $1\times 1$ Bernoulli entries with the appropriate optimal $p$ for 1's (as found above). By adding these up we again obtain a density profile, this time for the $1\times 1$ optimal Bernoulli measure at the resolution level of the block size $3\times 3$. The $3\times 3$ blocks can of course be generated using the optimal $2\times 2$ blocks as well and yet another density profile results. These three discrete plots are rendered as curves in Figure 5.

Perhaps the most notable feature here is the flattening of the distributions, as the block size increases i.e. the total block probabilities move towards the tails (while their means stay constant around $0.22$). The curves cross between density levels $1/3-4/9$: below this cross over the shorter range Bernoulli measures favor light $3\times 3$ blocks, above it they discount heavier blocks in comparison to the optimal $3\times 3$ Bernoulli measure.

When examined closer one will see that the total probability of $3\times 3$ blocks at a given density level essentially comes from at most three different kinds of local configurations (up to reductions above that is). These seem to be \lq\lq grown\rq\rq : when moving from density level $d$ to level $d+1/9$ the high probability blocks are generated by adding a (contiguous) 1 into an existing high probability block. This mechanism cannot prevail when the $3\times 3$ blocks are generated independently from smaller blocks. Consequently the small block curves in Figure 5. have suppressed tails. We expect this phenomenon to prevail in the higher order Bernoulli blocks as well and thereby to be a significant feature in the long range correlations of the measure of maximal entropy.

\vskip .4truein
\centerline{\hbox{
 \psfig{figure=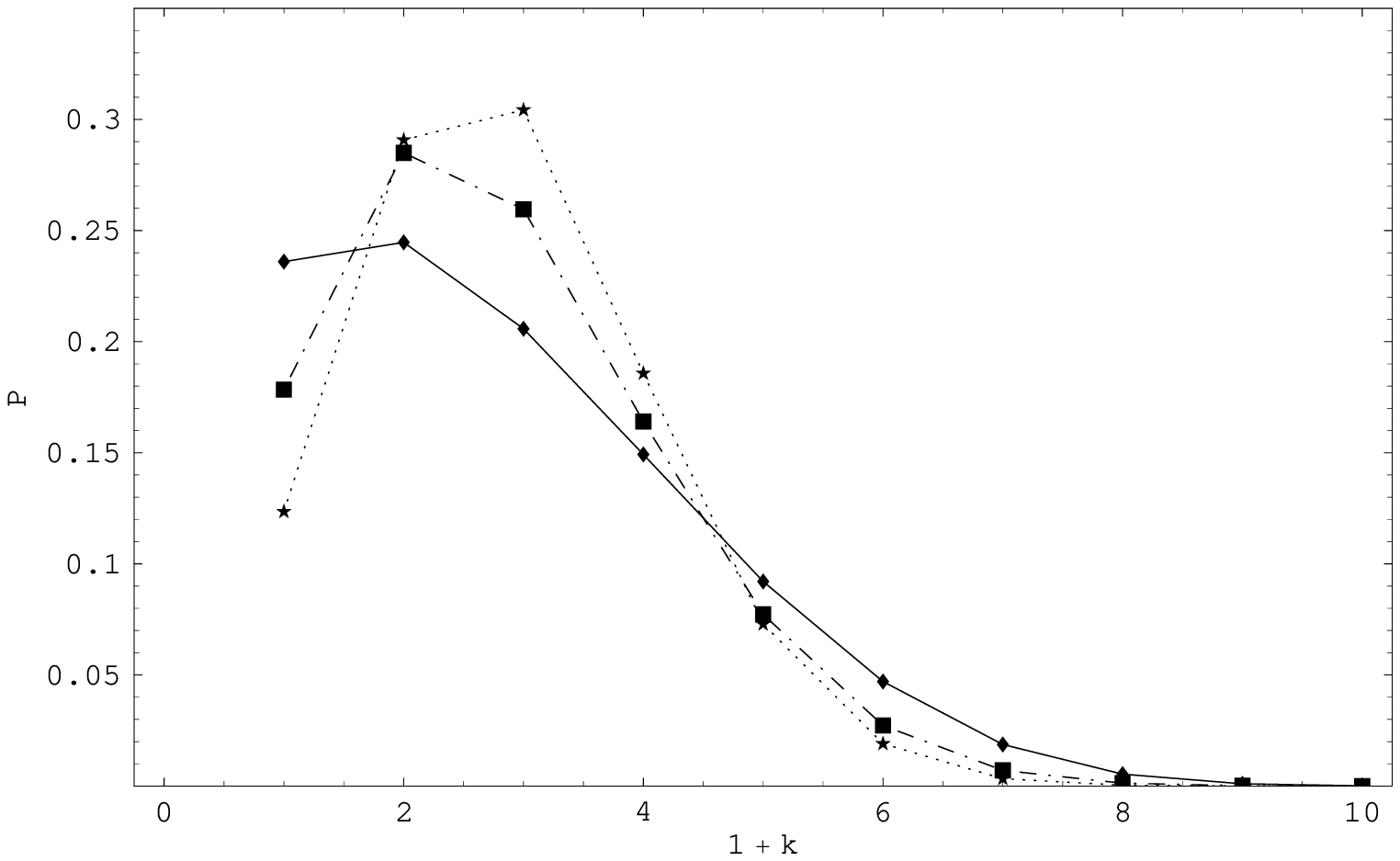,height=2in}
}}
\vskip .2truein
\nobreak
\noindent {\bf Figure 5.} $3\times 3$ -block occupation probabilities from Bernoulli blocks of size $3\times 3$ (diamond), $2\times 2$ (square) and $1\times 1$ (star). $k\in \{0,1,\ldots,9\}$ is the number of 1's in the block.

\vskip .6truein
\noindent {\subtitle References}
\vskip .3truein

\item{[B1]} Baxter, R.J.: {\sl Exactly solvable models in statistical mechanics}, Academic Press, 1982.

\item{[B2]} Baxter, R.J.: Planar lattice gasses with nearest-neighbor exclusion, {\sl Annals of Combinatorics}, {\bf 3}, no. 2-4, pp. 191-203, 1999, MR1772345 (2001h:82010).

\item{[DS]} Dobrushin, R., Shlosman, S., in {\sl Statistical Physics and Dynamical Systems}, J. Fritz, A. Jaffe and D. Szasz, eds., Birkh\"auser, pp. 347-370, 1985.

\item{[E]} Eloranta, K.: Dense packing on uniform lattices, {\sl J. of Stat. Phys.}, {\bf 130}, no. 4, pp. 741-55, 2008, {\tt arXiv:math-ph/0907.4247}, MR2387564 (2009a:52016).

\item{[MSS]} Milosevic, S., Stosic, B., Stosic, T.: Towards finding
exact residual entropies of the Ising ferromagnets, {\sl Physica A},
{\bf 157}, pp. 899-906, 1989.

\item{[RS]} Radulescu, D., Styer, D.: The Dobrushin-Shlosman Phase Uniqueness Criterion and Applications to Hard Squares, {\sl J. of Stat. Phys.}, {\bf 49}, pp. 281-95, 1987, MR0923859 (89d:82010).

\item{[W]} Walters, P.: {\sl Ergodic theory}, Springer, 1982.

\vfill
\eject
\end